\documentclass{article}
\usepackage{amssymb}
\newcommand{\eq}{\begin{equation}}
\newcommand{\en}{\end{equation}}

\newcommand{\rank}{\rm RANK}

\newcommand{\giv}{\,|\,}
\newcommand{\rev}{^\leftarrow}
\newcommand{\dec}{^\downarrow}

\newcommand{\MBR}{self-similar Markov}

\newcommand{\re}[1]{\mbox{(\ref{#1})}}
\newcommand{\rem}[1]{\mbox{\rm (\ref{#1})}}
\newcommand{\lb}[1]{\label{#1}}

\newcommand{\prob}{\mathbb P}

\newcommand{\ed}{ \stackrel{d}{=}}

\newcommand{\drift}{{\tt d}}
\def\endpf{\hfill $\Box$ \vskip0.5cm}

\newtheorem{theorem}{\large Theorem}
\newtheorem{proposition}[theorem] {\large Proposition}
\newtheorem{definition}[theorem]{\large Definition}
\newtheorem{corollary}[theorem]{\large Corollary}
\newtheorem{lemma}[theorem]{\large Lemma}

\begin{document}

\title{
Self-similar and Markov composition structures
\thanks{Research supported in part by N.S.F. Grant DMS-0405779
}
}
\author{Alexander Gnedin\thanks{Utrecht University; e-mail gnedin@math.uu.nl}
\hspace{.2cm}
and 
\hspace{.2cm}
Jim Pitman\thanks{University of California, Berkeley; e-mail pitman@stat.Berkeley.EDU} 
\\
\\
\\
\\
}

\date{\today}
\maketitle

\vskip0.5cm

\noindent
{\bf Abstract} 
The bijection between composition structures and random closed
subsets of the unit interval implies that the composition structures 
associated with $S \cap [0,1]$ for a self-similar random set 
$S\subset {\mathbb R}_+$ are those which are consistent 
with respect to a simple truncation operation. Using the standard coding
of compositions by finite strings of binary digits starting with a $1$, 
the random composition of $n$ is defined by the first $n$ terms of a
random binary sequence of infinite length.
The locations of $1$s in the sequence are the places visited by
an increasing time-homogeneous Markov chain on the positive integers if and
only if $S = \exp( -W )$ for some stationary regenerative random subset $W$
of the real line.
Complementing our study in previous papers, we identify 
self-similar Markovian composition structures associated with the 
two-parameter family of partition structures.

\section{Introduction}\label{intro}

A {\em composition of} $n$ is a sequence $\lambda=(\lambda_1,\ldots,\lambda_\ell)$ of some number $\ell$ of 
positive integer {\em parts} $\lambda_i$ with 
$\sum_{i = 1}^ \ell \lambda_i = n$.
We may regard $\lambda$ as a distribution of $n$ identical
balls in a row of $\ell$ boxes, with $\lambda_i$ the number of balls
in the $i$th box from the left end of the row.
Thus the composition $(2,4,1,2)$ of $9$ may be represented 
in {\em balls-in-boxes notation} as
\eq
\label{balls_boxes}
(2,4,1,2) \leftrightarrow [00]~[0000]~[0]~[00]
\en
or recoded in {\em binary notation} by replacing each ``$[0$" in the
balls-in-boxes 
notation by $1$ and ignoring each ``$]$", to obtain in this example 
$$
(2,4,1,2) \leftrightarrow 101000110.
$$
In general, the first digit in the binary notation of a composition
must be a $1$, but the remaining digits can be chosen freely,
so there are $2^{n-1}$ different compositions of $n$.
Two other notations will be useful. We write $\lambda\rev$ for the
{\em reversal} of $\lambda$ and $\lambda\dec$ for the 
{\em decreasing rearrangement of } $\lambda$, also called the
{\em partition derived from $\lambda$}. For instance
$$
(2,4,1,2) \rev = (2,1,4,2) \leftrightarrow [00]~[0]~[0000]~[00] \leftrightarrow 101100010
$$
$$
(2,4,1,2) \dec = (4,2,2,1) \leftrightarrow [0000]~[00]~[00]~[0] \leftrightarrow 100010101
$$
A {\em random composition of} $n$ is a random variable ${\cal C}_n$ with values in the set of all compositions of $n$.
We are interested in sequences of random compositions 
$({\cal C}_n)$ which are consistent as $n$ varies with respect to various {\em reduction operations}. In the balls-in-boxes description,
let the places of the $n$ balls be indexed from left to right by the set $[n]:= \{1, \ldots, n\}$. Let $(Y_n)$ be a sequence of random variables with $Y_n \in [n]$
for each $n$, with $Y_n$ independent of ${\cal C}_n$.
Let ${\cal C}_n^{-}$ be the composition of $n-1$ obtained by deleting
the ball in place $Y_n$ from the balls-in-boxes representation of ${\cal C}_n$.
For instance, if 
$$
{\cal C}_9 =(2,4,1,2) \leftrightarrow [00]~[0000]~[0]~[00] \leftrightarrow 101000110
$$
as above, and $Y_{9} = 6$, then 
$$
{\cal C}_9^{-} = (2,3,1,2) \leftrightarrow [00]~[000]~[0]~[00] \leftrightarrow 10100110.
$$ 
Whereas if instead $Y_{9} = 7$, then 
$$
{\cal C}_9^{-} = (2,4,2) \leftrightarrow [00]~[0000]~[00] \leftrightarrow 10100010.
$$
We say that the sequence of compositions $({\cal C}_n)$ 
is {\em $(Y_n)$-consistent} if there is
the equality in distribution
\eq
\lb{eqd}
{\cal C}_n^- \ed {\cal C}_{n-1} \mbox{ for every } n = 2,3, \ldots
\en
where ${\cal C}_n^{-}$ is ${\cal C}_n$ reduced by deletion of the ball in place
$Y_n$. 
Then, by Kolmogorov's extension theorem,
the sequence $({\cal C}_n)$ can be realised jointly
with $(Y_n)$ on a common probability space, so that the equality
in \re{eqd} also holds almost surely.
We say that such a realisation of $({\cal C}_n)$ is {\em strong $(Y_n)$-consistent}.
We are particularly concerned with the operations of 
{\em uniform, left} and {\em right} reduction corresponding to
$Y_n$ with uniform distribution on $[n]$, to $Y_n \equiv 1$, and to
$Y_n \equiv n$. That is, removal of a ball picked uniformly at random,
or the left-most ball, or the right-most ball. So we may call
$({\cal C}_n)$ {\em uniform-, left-} or {\em right-consistent} as the case may be.
Note that $({\cal C}_n)$ is left-consistent iff $({\cal C}_n\rev)$ is 
right-consistent.

\par There is an obvious bijection between sequences of distributions of 
${\cal C}_n$ which are right-consistent and probability distributions of 
infinite binary sequences $(\xi_1,\xi_2, \ldots)$ with $\xi_1 = 1$:
a strong right-consistent realisation of ${\cal C}_n$ in binary notation is 
the truncation $(\xi_1,\xi_2, \ldots, \xi_n)$ of the infinite binary sequence.
An alternate representation is obtained by replacing $(\xi_i)$ by the random set of 
positive integers $\{i \ge 1: \xi_i = 1\}$. 
Thus right-consistent sequences of compositions may be identified with random subsets of positive integers 
which contain $1$.

A uniform-consistent sequence of random compositions $({\cal C}_n)$ is 
also called a {\em composition structure} \cite{gnedin97, rcs}.
The corresponding sequence of random partitions $({\cal C}_n\dec)$ is 
then a {\em partition structure} in the sense of 
Kingman (see \cite{csp} for a survey and background).
That is to say,
\eq
\lb{eqd1}
({\cal C}_n^{\downarrow - })\dec \ed {\cal C}_{n-1} \dec \mbox{ for every } n = 2,3, \ldots.
\en
where the left side is the decreasing rearrangement of a reduction of
${\cal C}_n\dec$ by removal of uniformly chosen random ball.
Kingman gave a representation of partition structures which Gnedin refined as follows:

\begin{theorem}\label{gnedin} 
{\em \cite{gnedin97}}
Let $({\cal C}_n)$ be a composition structure.
Then there exists a random closed subset $Z$ of $[0,1]$ such that
a strong uniform-consistent realisation of $({\cal C}_n)$ can be constructed
as follows: let $(U_i)$ be a sequence of independent uniform $[0,1]$ variables,
independent of $Z$, and let ${\cal C}_n$ be
the sequence of sizes of equivalence classes among $U_1,\ldots,U_n$, 
listed left to right, as these points are classified
by the random equivalence relation $\sim$ induced by
$U_i \sim U_j$ for $i \ne j$ iff $U_i$ and $U_j$ fall in the same
interval component of $[0,1] \backslash Z$.
\end{theorem}

\noindent
{\bf Remark.} For consistency with further considerations in this paper 
we include the point $1$ in $Z$ only if $1$ is not an isolated point in $Z$.
Thus, if the rightmost interval of $[0,1]\setminus Z$ exists, it is semiopen.
\vskip0.5cm

\par
We are most interested in the case when $Z$ is {\em light}, meaning that
the Lebesgue measure of $Z$ equals $0$ almost surely.
The collection of component intervals of $[0,1] \backslash Z$ then defines
a {\em random interval partition} of $[0,1]$, that is a collection of
open subintervals of $[0,1]$, the sum of whose lengths is $1$.
The collection of lengths of component intervals of $[0,1] \backslash Z$, suitably
indexed, is then a {\em random discrete distribution} as studied in \cite{RDD}.
We regard the composition structure $({\cal C}_n)$ as a combinatorial representation of
either $Z$ or its associated interval partition, just as the partition structure 
$({\cal C}_n \dec)$ may be regarded as a combinatorial representation of 
the unordered collection of interval lengths.

In a series of previous papers \cite{slow, gnedin97, sieve, tsf, rcs, rps, gpyI, gpyII, bernoulli, RDD}, 
we have studied the composition structures, 
partition structures, random interval partitions, and random discrete distributions,
corresponding to various random subsets of $[0,1]$ of particular interest.
Here we tie together some threads from these previous studies, to show
how various analytic properties of the random subset $Z$ of $[0,1]$, 
which are natural from the perspective of continuous parameter stochastic processes,
correspond to various combinatorial properties of the associated composition
structure $({\cal C}_n)$.

A random closed subset $S$ of ${\mathbb R}_+$ 
is called {\em self-similar} (or {\em scale-invariant}) if
\begin{equation}\label{sss}
S\ed c\,S \mbox{ for all } c>0.
\end{equation}
In Section \ref{sec.ss} we establish the following result, which generalizes 
a construction introduced in \cite{bernoulli} in the case discussed in Example 2 below.

\begin{theorem}\label{general} 
For a sequence of distributions of random compositions $({\cal C}_n)$ 
the following two conditions are equivalent:
\begin{itemize}
\item $({\cal C}_n)$ is both uniform-consistent and right-consistent.
\item $({\cal C}_n)$ can be derived by uniform sampling from $S \cap [0,1]$ for some self-similar random closed subset $S$ of ${\mathbb R}_+$ .
\end{itemize}
When these conditions hold, a strong right-consistent version of
$({\cal C}_n)$ can be constructed as follows:
independent of $S$, let $\epsilon_1<\epsilon_2<\ldots$ be the atoms of a homogeneous 
Poisson point process (henceforth PPP) on ${\mathbb R}_+$, and let the binary representation
of ${\cal C}_n$ be the first $n$ digits of $(\xi_i)$ defined by $\xi_1=1$ and for $j>1$ 
$$\xi_j=1([\epsilon_{j-1}\,,\,\epsilon_j]\cap S\neq\emptyset).$$
\end{theorem}

Note that there are two quite different 
realisations of $({\cal C}_n)$, one that is strong uniform-consistent, obtained by 
uniform sampling from $S \cap [0,1]$ as described in Theorem \ref{gnedin}, and one that is 
strong right-consistent, obtained by Poisson sampling from $S \subseteq [0, \infty[$ as 
described in Theorem \ref{general}. Obviously,
it impossible to construct $({\cal C}_n)$ to be simultaneously 
strong uniform-consistent and strong right-consistent.

\vskip0.5cm
\noindent
{\bf Example 1} \cite{abt, dj91, tsf, rcs} Let $(\xi_1,\xi_2,\ldots)$ be a random Bernoulli string with independent digits 
and distribution
$${\mathbb P}(\xi_j=1)=1-{\mathbb P}(\xi_j=0)=\theta/(j+\theta-1),$$
where $0<\theta<\infty$ is a parameter. 
Let ${\cal C}_n$ be encoded by the first $n$ digits, so 
$({\cal C}_n)$ is strong right-consistent by construction.
It is elementary that 
${\cal C}_n$ is 
distributed according to the
formula
\eq\label{pEw}
{\mathbb P}({\cal C}_n=(\lambda_1,\ldots,\lambda_\ell))={\theta^\ell n!\over (\theta)_n}\prod_{j=1}^\ell
{1\over\Lambda_j}
\en
where $(\theta)_n=\theta(\theta+1)\cdots(\theta+n-1)$ and $\Lambda_j=\lambda_1+\cdots+\lambda_j$.
This is a variant of the Ewens sampling formula \cite{dj91, tsf},
which gives the distribution of the sizes of blocks, in reverse 
size-biased order, of a Ewens partition of $n$ with parameter $\theta$.
In the case $\theta=1$ the sequence $(\xi_i)$ results from encoding 
the cycle partition of a uniform random permutation of $[n]$ by
Feller coupling \cite{abt}. 
This sequence also appears in the theory of extremes
as the sequence of record indicators of independent identically distributed observations 
with continuous distribution \cite{Bunge}.
The uniform-consistency of $({\cal C}_n)$ was observed in \cite{dj91}.
It is known that the random set $Z$ in Kingman's representation is 
the restriction to $[0,1]$ of the self-similar random set
$S$ which is the union of $\{0\}$ and the set of points of a scale-invariant Poisson process 
on $[0,\infty[$, with intensity $\theta \, dx/x\,,\,x > 0 $.
Properties of this scale-invariant Poisson process are reviewed in \cite{arratia}.
Two trivial composition structures 
appear as limiting cases for $\theta \downarrow 0$ and 
$\theta \uparrow \infty$.

\vskip0.5cm
\noindent
{\bf Example 2} \cite{bernoulli, rcs} Let $\xi_j=1(R_k=j~{\rm for~some~}k \ge 0)$ where $R_0 = 1$ and 
$R_k=1+X_1+\cdots+X_k$ is the discrete 
renewal process derived from independent and identically
distributed $X_j$ with
$${\mathbb P}(X_j =r)=(-1)^{r-1}{\alpha\choose r},~~r=1,2,\ldots\,$$
where $0< \alpha< 1$. The corresponding right-consistent sequence
of compositions has distribution
\eq\label{p-ren}
{\mathbb P}({\cal C}_n=\lambda)=\lambda_\ell \,\alpha^{\ell-1}\prod_{j=1}^\ell {(1-\alpha)_{\lambda_j-1}\over\lambda_j!}\,.
\en
That this sequence of compositions $({\cal C}_n)$ is both right-consistent and
uniform-consistent was shown in \cite{bernoulli}, where the two different 
strong consistent constructions of Theorems \ref{general} and \ref{gnedin} were given for 
$S$ the self-similar zero set of a Bessel process of dimension $2-2\alpha$, and
$Z = S \cap [0,1]$.
The trivial cases appear again as limits for $\alpha=0{\rm ~or~}1$.
\vskip0.5cm

\par It was shown by J. Young \cite{young} that 
no other choice of distribution either for a sequence of independent
Bernoulli variables (as in Example 1), or for a renewal sequence with independent
spacing between $1$'s (as in Example 2) yields a right-consistent sequence of 
compositions $({\cal C}_n)$ such that $({\cal C}_n\dec)$ is a partition structure.
As the latter condition is weaker than uniform-consistency of $({\cal C}_n)$,
no more right-consistent composition structures can be obtained from these constructions.
However, we will show that a construction adopted from 
\cite{RDD, young, rcs} allows an interesting extrapolation of the above examples 
to obtain a right-consistent composition structure with two parameters 
$(\alpha, \theta)$ with $0 \le \alpha<1$ and $\theta>-\alpha$, corresponding to the
two-parameter Ewens-Pitman family of partition structures.

\par
To emphasise the general correspondence between composition structures and
random sets provided by Theorem \ref{gnedin}, we use terminology for composition
structures to reflect properties of their associated random sets. 
So we prefer the term {\em self-similar} rather than right-consistent
for the composition structure $({\cal C}_n)$ obtained by uniform sampling from $S \cap [0,1]$ for a self-similar random set $S$. 
In \cite{rcs} we described the random sets associated
with composition structures $({\cal C}_n)$ with the following
{\it left-regenerative} property: for every $n$ and $1 \le x \le n$,
conditionally given the leftmost part of ${\cal C}_n$ is $x$, 
the remaining composition of $n-x$ is a distributional copy of ${\cal C}_{n-x}$.
Here we find it more convenient to work with the {\it right-regenerative} property, 
defined in the same way with the rightmost part instead of the leftmost part. 
Evidently, $({\cal C}_n)$ derived by uniform sampling from $Z$ is right-regenerative 
iff $({\cal C}_n\rev)$ derived by uniform sampling 
from $1-Z$ is left-regenerative. So the main result of \cite{rcs} can be 
restated as follows: a composition structure $({\cal C}_n)$ is right-regenerative 
iff $({\cal C}_n)$ is derived by uniform sampling from $e^{-W}$ for $W$ 
a regenerative random subset of $[0,\infty[\,$.
Here we distinguish a class of {\em Markov} composition structures such that the
binary representation of ${\cal C}_n$ has $1$'s at the places visited by a decreasing
Markov chain on $[n]$ with some transition matrix $q$ which does not depend on
$n$, and some initial distribution $q_*(n,\cdot)$ on $[n]$. These turn out to
be derived by uniform sampling from $e^{-W}$ for $W$ a delayed regenerative random 
subset of $[0,\infty[\,$. In the special case when $W$ is a {\it stationary} regenerative set,
$e^{-W}$ is the restriction to $[0,1]$ of a self-similar random subset of $[0,\infty[\,$.
The {\em self-similar Markov} compositions so obtained 
turn out to be those whose infinite binary representation has $1$'s at the places visited 
by an increasing Markov chain on ${\mathbb N}$. 
Finally, extending our study in \cite{rcs}, we introduce self-similar Markov composition 
structures associated with the two-parameter Ewens-Pitman family of partition structures.

\section{Self-similar composition structures}
\label{sec.ss}

\subsection{Composition probability function}

The distribution of a random composition ${\cal C}_n$ of integer $n$ is a nonnegative function 
$$p(\lambda)={\mathbb P}({\cal C}_n=\lambda)$$
on compositions $\lambda$ of $n$
which satisfies $\sum_{\lambda: |\lambda|=n}p(\lambda)=1$.
Here and henceforth $|\lambda|$ denotes the sum of parts of a composition $\lambda$.
For a general sequence of random compositions $({\cal C}_n)$ these marginal distributions are
described by a
{\it composition probability function} (CPF)
defined for all compositions of integers. 
A sequence of compositions is $(Y_n)$-consistent iff the CPF 
satisfies
a linear recurrence of the form
\begin{equation}\label{S}
p(\lambda)=\sum_{\mu: |\mu|=|\lambda| +1}
p(\mu) \, \varkappa(\mu, \lambda)
\end{equation}
where $\varkappa(\mu,\lambda)$ for $\mu$ with $|\mu| = n$ is
a matrix describing the transition probabilies from compositions of $n$ to 
compositions of $n-1$ determined in the balls-in-boxes representation by removal of a ball
from place $Y_n$.
See $\cite{gnedin97, rcs}$ for 
details in the case of uniform-consistency when $Y_n$ has uniform distribution on $[n]$.
For $({\cal C}_n)$ that is right-consistent, the recurrence is just
linear relation
\begin{equation}\label{T}
p(\lambda_1,\ldots,\lambda_\ell)=p(\lambda_1,\ldots,\lambda_\ell+1)+p(\lambda_1,\ldots,\lambda_\ell,1).
\end{equation}

\subsection{Proof of Theorem 2}

We start by remarking that
the distribution of a self-similar set $S$ is uniquely determined by the distribution of 
its restriction $Z$ to $[0,1]$, 
which satisfies 
the condition equivalent to (\ref{sss}):
\eq
\label{sss1}
c \,(S\cap [0,1])\ed S\cap[0,c]\,,~~~{\rm for~~}\, 0<c<1\,. 
\en
This follows from the known fact that the distribution of a {\it stationary} set $W\subset {\mathbb R}$
(invariant under shifts) is determined by the distribution of $W\cap {\mathbb R}_+$,
and we can transform a self-similar $S$ into a stationary set
$W:=-\log S$.

For $n=1,2,\ldots$ let 
$Z_n\subset[0,1[$ be a finite set
encoding ${\cal C}_n$ via the correspondence 
$(\lambda_1,\ldots,\lambda_\ell)\to \{0, \Lambda_1/n,\ldots,\Lambda_{\ell-1}/n\}$ 
for $\Lambda_j=\lambda_1+\ldots+\lambda_j$.
Assuming now that we are working with a 
strong uniform-consistent realisation of $({\cal C}_n)$,
by the law of large numbers \cite{gnedin97}
the Hausdorff distance between $Z_n\cup \{1\}$
and $Z\cup \{1\}$ goes to $0$ with probability $1$.
The Hausdorff distance between $Z_n$ and $Z$ also goes to $0$. 
This can be shown by considering the last block of the composition, which 
has a positive frequency if and only if $1$ is not an accumulation
point for $Z$.
In this sense,
$Z_n\to Z$ a.s., hence also 
$Z_{\lfloor nx\rfloor}\to Z$ a.s.
for every $x\in \,]0,1[\,$.
Translating the 
truncation property in terms of $Z_n$'s we obtain
$$
Z_n\cap \left[0\,,\,{{\lfloor nx\rfloor}\over n}\right[\,\ed {{\lfloor nx\rfloor}\over n}Z_{\lfloor nx\rfloor}.
$$
For $n\to\infty$ the left side converges to $Z\cap [0,x[$ a.s., while the right side converges 
to $x\,Z$ a.s., hence the limits must have the same distribution. 
This means that $Z$ is self-similar, by \re{sss1}.
The strong right-consistent represention is obtained by noting that the scaling
$(S, \epsilon_1,\ldots,\epsilon_{n})\to (S/\epsilon_{n+1}, \epsilon_1/\epsilon_{n+1},\ldots,\epsilon_n/\epsilon_{n+1})$
transforms $S$ to a copy of itself (by self-similarity) and maps 
the first $n$ Poisson points to the increasing sequence of $n$ uniform order 
statistics. 
\endpf

\subsection{Some definitions}
\par We call $Z$ {\it heavy} if $Z$ has positive Lebesgue measure with nonzero probability, and we call $Z$
{\it light} otherwise.
The set $Z$ can be discrete (as in Example 1) or perfect (as in Example 2) or neither discrete nor perfect.
\par For $x\in {\mathbb R}_+$ introduce 
$$G_x :=\sup(Z\cap[0,x])\,,~~A_x :=x-G_x\,,~~
D_x := \inf(Z\cap \,]x,\infty[).$$
The {\it age process} $(A_x\,,\,x\geq 0)$ uniquely determines $Z$.
In the event $x\in Z$ we have $G_x=x$ and $A_x=0$, while in the event $x\notin Z$ the point $x$ is covered
by an open {\it gap} $]G_x, D_x[\,\subset {\mathbb R}_+\setminus Z$.
The interval $]G_1,1[$ of length $A_1$ is called the {\it meander}.
If $Z$ is heavy the meander may be empty with positive probability (then $A_1=0$), 
while for light $Z$ the meander is nondegenerate (and $A_1>0$ a.s.).

\section{Block counts, meander and the tagged interval}
\label{block}

\subsection{The structural distribution}

For $Z\subset [0,1]$ a random closed set and $U$ a uniform random point independent of $Z$ let 
$V$ be the size of gap in $Z$ covering $U$ in case $U\in [0,1]\setminus Z$, and let $V=0$ in case
$U\in Z$. The gap covering a random point is sometimes called the {\it tagged} interval.

\par Let $(V_j)$ be the decreasing sequence of lengths of gaps comprising $[0,1]\setminus Z$,
so that $\sum_j V_j\leq 1$ and $1-\sum_j V_j$ is the Lebesgue measure of $Z$. 
If all the $V_j$'s are pairwise distinct with probability one, then 
$${\mathbb P}(V=V_j\,|\,V_1,V_2,\ldots)=V_j\,,~~~{\mathbb P}(V=0\,|\,V_1,V_2,\ldots)=1-\sum_j V_j$$
So $V$ may be called a {\em size-biased pick} from the sequence of lengths.

\par Suppose now that the composition structure $({\cal C}_n)$ is derived by uniform
sampling from $Z$.
The distribution of $V$ is called the {\it structural distribution}
(of $({\cal C}_n)$ , or of the associated partitition structure, or of the
associated random discrete distribution of interval lengths).
Recall that $p$ denotes the CPF of $({\cal C}_n)$. Observe that
$$p(n) = {\mathbb E}\,V^{n-1}$$
because ${\cal C}_n$ equals the one-part composition $(n)$ when $n-1$ further uniform points
hit the interval of length $V$ found by $U$. Other relations of this type are
\begin{eqnarray}
\label{mur1}
\mu_{n,r}&:=&{\mathbb E}\,K_{n,r}={n\choose r}{\mathbb E}\left(V^{r-1}(1-V)^{n-r}\right)\,,~~~~~1\leq r\leq n.\\
\label{mur}
\mu_n&:=&{\mathbb E}\,K_{n}=\sum_r \mu_{n,r}=
{\mathbb E}\left({1-(1-V)^n\over V}\,1(V>0)\right)+n\,{\mathbb P}(V=0)
\end{eqnarray}
where $K_{n,r}$ is the number of parts of ${\cal C}_n$ or size $r$,
and $K_n=\sum_r K_{n,r}$ is the number of parts of 
${\cal C}_n$.
Note that $(K_{n,r},1 \le r \le n)$ is a standard encoding of 
${\cal C}_n\dec$,
the random partition of $n$ induced by ${\cal C}_n$.
\begin{theorem}\label{A=V}
{\rm \cite{RDD}} 
Suppose that $S\subset {\mathbb R}_+$ is self-similar, and let
$Z:= S \cap [0,1]$.
Let $A_1$ be the length of the meander interval of $Z$, that is
the rightmost gap in $[0,1] \backslash Z$, with $A_1 = 0$ if $1\in Z$,
and let $V$ be the length of the component interval of $[0,1] \backslash Z$ which 
contains $U$ independent of $Z$, with $V = 0$ if $U\in Z$.
Then $A_1$ has the same distribution as $V$.
\end{theorem}
{\it Proof.} Let $\rho_n$ be the rightmost in the sample of $n$ uniform points.
Given $A_{\rho_n}$, with probability
$(A_{\rho_n}/\rho_n)^{n-1}$ the remaining $n-1$ sample points fall in the
same gap of $Z$ as $\rho_n$.
By self-similarity, $A_{\rho_n}/\rho_n\ed A_1$. So
the probability that all $n$ points fall in the same gap is 
$$p(n)={\mathbb E}\left(A_{w}/w\right)^{n-1}={\mathbb E}A_1^{n-1}.$$
Comparing with $p(n)={\mathbb E}V^{n-1}$ we arrive at the conclusion, since 
a probability distribution on $[0,1]$ is determined by its moments.
\endpf
\noindent 
Note that the event $(A_1=V)$ has probability ${\mathbb E}A_1={\mathbb E}\,V=p(2)$.
Pitman and Yor \cite{RDD} went further to distinguish 
a {\it strong sampling property} for the meander 
\eq\label{ssp}
{\mathbb P}(A_1=V_k\,|\,V_1,V_2,\ldots)=V_k\,,~k=1,2,\ldots
\en
meaning the condition that the meander length is a size-biased pick from all lengths.
This property holds in some cases 
({\it e.g.} for $Z$ in Examples 1 and 2, and in the setup of Theorem \ref{thmalth} below) but 
does not hold in general.

\subsection{The last part and the tagged part of composition}

Theorem \ref{A=V} implies that for self-similar composition structure, 
as $n\to\infty$,
the frequency of the last block of ${\cal C}_n$
has approximately the same distribution as the frequency of the block selected by a size-biased pick.
A stronger fact is true: a similar
identity holds for each $n$, and not only asymptotically.
This was already observed in \cite[Proposition 11 (i)]{bernoulli}
in the case of renewal strings in Example 2.
Intuitively, since both uniform- and right- reduction transform ${\cal C}_n$ into 
a composition with the same distribution, 
it is natural to expect that the sizes of reduced parts
have the same distribution.

\begin{theorem}\label{S=L} 
For a composition structure $({\cal C}_n)$,
let $P_n$ denote the size of a random part of ${\cal C}_n$
which given ${\cal C}_n$ is selected with probability proportional to size.
Let $L_n$ be the size of the last part of ${\cal C}_n$.
If $({\cal C}_n)$ is right-consistent then $P_n\ed L_n\,$ for all $n$.
\end{theorem}

This follows immediately from the following Lemma.

\begin{lemma}\label{pi} 
Let ${\cal C}_{n}$ and ${\cal C}_{n-1}$ be two random compositions of
$n$ and $n-1$ respectively, defined on a common probability space
in such a way that ${\cal C}_{n-1}$ is obtained from ${\cal C}_{n}$ by 
removal of a single ball in the balls-in-boxes representation.
Let $(\omega_{n,r}, 1\leq r\leq n)$ be the distribution of the number of balls
in the same box of ${\cal C}_{n}$ as the ball removed,
and let $\mu_{n,r}$ and $\mu_{n-1,r}$ be the expected numbers of boxes containing
$r$ balls for ${\cal C}_n$ and ${\cal C}_{n-1}$, respectively, as above in {\rm (\ref{mur1})}.
Then the distribution $(\omega_{n,r}, 1\leq r\leq n)$ is determined by the
distributions of the partitions generated by ${\cal C}_{n-1}$ and ${\cal C}_{n}$ according
to the formulas $\omega_{n,n}=\mu_{n,n}$ and
$$
\omega_{n,r}-\omega_{n,r+1}=\mu_{n,r}-\mu_{n-1,r}
~~~ ( 1\leq r\leq n-1 ) .
$$
\end{lemma}
{\it Proof.} Follow the evolution of $K_{n,r}$ as $n$ varies. This variable
increases by $1$ when a ball is chosen in a box of $r+1$ balls 
(which is impossible for $r=n$), and decreases by $1$ when a ball is chosen in a box of $r$ balls.
The probabilities of these events are $\omega_{n,r+1}$ and $\omega_{n,r}$, respectively. 
In all other cases the sampling does not affect $K_{n,r}$. 
The formula for expected increments follows.
\endpf 
\noindent 
\noindent
Note that Theorem \ref{A=V} follows from Theorem \ref{S=L}
by the law of large numbers. 
A discrete analogue of (\ref{ssp}) holds for compositions in Examples 1 and 2: conditionally
given ${\cal C}_n \dec$, the last
part $L_n$ of ${\cal C}_n$ is a size-biased pick from all parts.

\subsection{A characterisation of structural distributions}

The following characterisation of structural distributions
is a minor extension of \cite[Condition 1]{RDD} to include the heavy case.

\begin{theorem}\label{strss} The structural distribution of the interval partition
derived from a self-similar random set $S$ has the form
\begin{equation}\label{form}
{\mathbb P}(V\in {\rm d}x)={\widetilde{\nu}[x,1]\over (\drift+{\tt m})\,(1-x)}\,{\rm d}x + {\drift\over\drift+{\tt m}}\,
\delta_0({\rm d}x)
\end{equation}
where 
$\widetilde{\nu}$ is a measure on $]0,1]$ satisfying
$$
{\tt m}:=\int_0^1 |\log(1-x)|\widetilde{\nu}({\rm d}x)<\infty\,
$$
and $\drift$ is a nonnegative constant.
Thus, the structural distribution may have an atom at $0$, and otherwise has a density $\phi(x),\,0<x\leq 1$,
such that $(1-x)\phi(x)$ is decreasing.
The data $(\drift, \widetilde{\nu})$ are determined uniquely up to a positive factor.
\end{theorem}
{\it Proof.} Assume the normalisation ${\tt m}=1$.
Let $W=-\log S$ be stationary and $X$ be a random variable whose distribution coincides with
the conditional distribution of the size of the gap of $W$ covering $0$ 
given this size is positive. By the ergodic theorem, the part of the gap on the positive halfline
is distributed like $X U$, with $U$ uniform $[0,1]$ independent of $X$. 
The conditional distribution of $A_1$ given $A_1>0$ is then the same as for $e^{-X U}$,
which implies along the lines of the argument in \cite[Section 4]{RDD} that
$A_1$ has a density written as $\widetilde{\nu}[x,1]/(1-x)$. The unconditional distribution in the form (\ref{form}) 
follows by defining $\drift$ from 
$${\mathbb P}(A_1=0)={\mathbb P}(0\in W)={\drift\over 1+\drift}$$
(where the middle term is the the long-run Lebesgue measure of $W$ per unit length).
\endpf

\par The structural distribution also accounts for 
some functionals of self-similar composition structures which involve the 
ordering of parts. For a self-similar composition structure
$({\cal C}_n)$ with binary representation $(\xi_1\xi_2\ldots)$,
define the {\it potential function} 
$$g(j)={\mathbb P}(\xi_j=1).$$
In terms of balls-in-boxes, this is the probability, for each $n \ge j$, that 
the $j$th ball of ${\cal C}_n$ falls in a different box from its predecessor.
Note that for a composition structure which was not right-consistent, 
the analogous quantity would typically depend on $n$ as well as $j$.
In terms of $\mu_n:= {\mathbb E} K_n$ and the moments of the structural distribution, we read
from \re{mur} that
$$g(j)=\mu_j-\mu_{j-1}={\mathbb E}\,(1-V)^{j-1}\,,~~~\mu_n =\sum_{j=1}^n g(j)$$ 
(where $\mu_0=0$).

\vskip0.5cm
\noindent
Nacu \cite{nacu} proved that the probability law
of the general exchangeable partition of ${\mathbb N}$ is uniquely determined by 
the distribution of the sequence $(I_j)$ of indicators of 
minimal elements of the blocks (so $I_7=1$ means that $7$ is the minimal element in some block).
In terms of Kingman's representation, $I_j=1$ each time $U_j$ discovers a new 
gap in $Z$ or hits $Z$.
For the Ewens composition structure of Example 1,
$(I_j)$ has the same distribution as $(\xi_j)$.
For a general self-similar composition structure, the sequences 
are differently distributed (as {\it e.g.} in Example 2), but
the right-consistency of ${\cal C}_n$ implies that
$${\mathbb P}(I_j=1)={\mathbb P}(\xi_j=1)=g(j)$$
because 
$$
I_1+\cdots+I_n = \xi_1+\cdots+\xi_n = K_n,
$$
the number of parts of ${\cal C}_n$.

\subsection{A fragmentation product}
\label{fp}
The following operation on self-similar sets 
generalises the one found in \cite{rcs, RDD, young}.
Let $Z\subset{\mathbb R}_+$ be self-similar and
independent of $Z$. Let $(M_{j})$ be independent copies of the same random closed 
set $M\subset[0,1]$.
For each gap in $Z$ with left-point $z_j\in Z$ and size $s_j$ fit the set $z_j+s_jM_j$ in this gap,
and take the union of $Z$ and all these scaled shifted copies of $M$. 
Then the result $Z\otimes M$ (read {\it $Z$ fragmented by $M$}) is easily shown to
be self-similar.
For example, when $M=\{1/2\}$ the set $Z\otimes M$ is obtained by adding the midpoint
for each gap in $Z$.

\par The operation has an analogue 
in terms of composition structures (as in \cite{young}). 
For two composition structures 
$({\cal C}_n)$ and $({\cal C}_n')$, for each $n$, 
break 
the generic part of ${\cal C}_n$, say $r$, into smaller parts
according to an independent copy of ${\cal C}_r'$.
The resulting sequence of compositions is a right-consistent composition structure provided
$({\cal C}_n)$ is so,
and the corresponding self-similar random set is the fragmentation product $Z\otimes M$ of
the sets in Kingman's representation of $({\cal C}_n)$ and $({\cal C}_n')$.

\section{Markovian composition structures}
\label{mcs}

\subsection{Decrement matrices}
The following extension of the concept of a {\it regenerative
composition structure} introduced in \cite{rcs} 
extends our study in that paper and prepares
for the results in the next section.

\begin{definition}
\label{Mcs}
{\rm A composition structure $({\cal C}_n)$ is called {\it Markovian} 
if for some infinite transition probability matrices 
$$(q(n:m),\,1\leq m\leq n<\infty)\quad {\rm and}\quad
(q_*(n:m),\,1\leq m\leq n<\infty)$$
the distribution of each ${\cal C}_n$ is given by the product formula
\begin{equation}\label{prod}
p(\lambda)=
q_*(n:\lambda_\ell)\prod_{k=1}^{\ell-1}q(\Lambda_{k}\,:\,\lambda_{k})\,,~~~~
\end{equation}
where $\lambda=(\lambda_1,\ldots,\lambda_\ell)$ is a composition of $n$, 
and $\Lambda_k=\lambda_1+\ldots+\lambda_k$ for $k\leq\ell$.}
\end{definition}

\par Formula (\ref{prod}) has the following interpretation. Imagine a {\it decreasing} time-homogeneous Markov chain 
$Q_n^\downarrow=(Q_{n,t}^\downarrow,\,t=0,1,\ldots)$ with state-space $\{1,2,\ldots,n\}$ and terminal absorbing state $1$.
The chain has initial distribution
$${\mathbb P}(Q_{n,0}^\downarrow=j)=q_*(n:n-j+1)\,,\quad j=1,\ldots,n$$
and it jumps from state $j$ ($2\leq j\leq n$) to $i$ ($1\leq i<j$) with probability $q(j-1:j-i)$. 
We call $q_*$ and $q$ {\it decrement matrices}.
In these terms, a random composition of $n$ can be identified with a path of $Q_n^\downarrow$, {\it i.e.}
the binary representation of ${\cal C}_n$ (for fixed $n$)
is obtained by writing $1$'s in the positions visited by $Q_n^\downarrow$.

\par In the case $q=q_*$ the formula (\ref{prod}) defines a regenerative composition structure, as introduced
in \cite{rcs}.
As mentioned in the Introduction, to fit in the present framework, 
the convention in that paper regarding the ordering of blocks should be reversed.

\begin{lemma}\label{mark-q}
For a Markovian composition structure we have 
\begin{eqnarray}\label{q}
q(n:r)={r+1\over n+1}\,q(n+1:r+1)+{n+1-r\over n+1}\,q(n+1:r)+{1\over n+1} \,q(n+1:1)\,q(n:r)\\
\label{q*}
q_*(n:r)={r+1\over n+1}\,q_*(n+1:r+1)+{n+1-r\over n+1}\,q_*(n+1:r)+{1\over n+1}\, q_*(n+1:1)\,q(n:r).
\end{eqnarray}
Conversely, if two nonnegative matrices $q_*$ and $q$ satisfy these recursions and 
$q_*(1:1)=q(1:1)=1$
then they define a Markovian composition structure via {\rm (\ref{prod})}.
\end{lemma}

\noindent
{\it Proof.} The recursions follow from (\ref{prod}) and uniform consistency. When a composition is reduced by sampling
the last block either remains unaltered or, in the case the last block is a singleton
and gets deleted, coincides with the second-last block, whence (\ref{q*}).

\par The first recursion is familiar from \cite{rcs, rps}, but proving it under the more general assumption
(\ref{prod})
requires more algebra. For $(a,b,c)$ a composition of $n$ use uniform consistency to obtain
\begin{eqnarray*}
p(a,b,c)= {1\over n+1} p(1,a,b,c)+ {a+1\over n+1}p(a+1,b,c)+
{1\over n+1}p(a,1,b,c)\\
+{b+1\over n+1}p(a,b+1,c)+{1\over n+1}p(a,b,1,c) +{c+1\over n+1}p(a,b,c+1)+{1\over n+1}p(a,b,c,1)\,.
\end{eqnarray*}
Group the first three terms in the right side as
$${a+1\over n+1} q_*(n+1:c)q(n+1-c:b)p(a)$$ 
and factor all other terms using (\ref{prod}). 
Factor the left side as
$$p(a,b,c)=q_*(n:c)q(a+b:b)p(a)$$
and express $q_*(n:c)$ through $q(n+1:\,\cdot\,)$ using (\ref{q*}). Cancelling common terms and factors
yields (\ref{q}).
\par The converse is checked as in \cite[Proposition 3.3]{rcs}.\endpf


\subsection{Kingman's representation}

Let $(Y_t, t \geq 0)$ be a {\it subordinator} (with $Y_0=0$),
meaning an increasing L{\'e}vy process. 
Let $0\leq X\leq \infty$ be a random variable, independent of $(Y_t)$ and satisfying ${\mathbb P}(X<\infty)>0$.
We call the process
$(X+Y_t,\,t\geq 0)$ a {\it delayed subordinator}, and call its closed range $W$ a {\it delayed regenerative set}.
The distribution of $W$ 
determines that of $X$ (because $X=\min W$) and determines the L{\'e}vy parameters $(\nu,\drift)$
up to a positive factor (since given $X<\infty$
the set $W-X$ is regenerative). 
Introduce the L{\'e}vy-Khintchine exponent 
\begin{equation}\label{Phis}
\Phi(s)=\drift s+\int_0^\infty (1-e^{-sy})\,\nu({\rm d}y)\,,~~~s\geq 0\,,
\end{equation}
its two-parameter extension
\begin{equation}\label{Phinm}
\Phi(n:m)={n\choose m}\sum_{j=0}^m (-1)^{j+1}{m\choose j} \Phi(n-m+j)\,,~~~~~1\leq m\leq n,
\end{equation}
and the moments
\begin{equation}\label{Psinm}
\Psi(n:m)={n\choose m}{\mathbb E}\left( (A_1^m (1-A_1)^{n-m}\right)\,,~~~~0\leq m\leq n
\end{equation}
where $A_1=1-\exp(-X)$.

\begin{theorem}\label{mark-gen}
A composition structure $({\cal C}_n)$ is Markovian if and only if it can be derived by
uniform sampling from $Z=\exp (-W)$, with $W$ being a delayed regenerative set. Explicitly, the distribution of $({\cal C}_n)$
is given by the product formula with decrement matrices
\begin{eqnarray}\label{qPhi}
q(n:m)&=&{\Phi(n:m)\over \Phi(n)}\\
\label{q*Psi}
q_*(n:m)&=&\Psi(n:0)q(n:m)+\Psi(n:m)\,.
\end{eqnarray}
\end{theorem}
{\it Proof.} The argument for the `if' part follows the same line as in \cite[Theorem 5.2 (i)]{rcs}. For the `only if' part 
let 
${\cal C}_n$ be derived by uniform sampling from
the random closed set $Z\subset [0,1]$.
Assume first that $G_1<1$ a.s. for $G_1=\sup Z\cap [0,1[\,$. 
Let $L_n$ be the last part of ${\cal C}_n$. Given $n-L_n=m$ let ${\cal C}_m'$ be a composition of
$m$ obtained by deleting the last part of ${\cal C}_n$. Note that this definition does not depend on $n>m$ and that
by Lemma \ref{mark-q} and \cite[Proposition 3.3]{rcs}, hence
$({\cal C}_m')$ is a regenerative composition structure.

\par Let $Z_m'$ be discrete random sets encoding ${\cal C}_m'$, 
as in the proof of Theorem \ref{general}, but with $1$ appended to $Z_m'$. 
The set $Z_n$ (not containing $1$) encoding 
${\cal C}_n$ can be represented as
$$Z_n\stackrel{d}{=}\left({n-L_n\over n}\right) Z_{n-L_n}'$$
where $L_n$ and $(Z_m')$ are independent. By the law of large numbers \cite{gnedin97}
$Z_n$ converge to $Z$, while
by \cite[Theorem 5.2 (ii)]{rcs} $Z_m'$ converge to some set $Z'=\exp(-W')$ with regenerative $W'\subset [0,\infty]$.
As $n\to\infty$ the law of the large numbers ensures that $1-L_n/n\to G_1$ a.s., hence in the limit we have
$$Z\stackrel{d}{=}G_1\,Z'$$
where $G_1$ and the set $Z'$ are independent. Hence the set $W=-\log Z$ is delayed regenerative.
\par The case ${\mathbb P}(G_1=1)>0$ is treated similarly. This can be viewed as a mixture (over $q_*$) of the trivial
one-block composition structure and another Markovian one.
\endpf

\section{Self-similar Markov composition structures}

\subsection{Markov sequences}
Let $Q^{\uparrow}=(Q^\uparrow_t, t=0,1,\ldots)$
be a time-homogeneous {\it increasing} Markov chain with the state-space 
$\{1,2,\ldots\}$ and the initial state $Q_0=1$. Define a string $\xi_1\xi_2\ldots$
by identifying the positions of $1$'s with the 
sequence of sites visited by $Q^\uparrow$:
$$\xi_j=1(Q_t=j~{\rm for~some~}t).$$
This defines a right-consistent sequence of compositions
$({\cal C}_n)$, so that each
${\cal C}_n$ encodes path of $Q^\uparrow$ killed before crossing level $n$.
We will consider such compositions which are also uniform-consistent, in which case 
(in view of Theorem \ref{general}) we will
call $({\cal C}_n)$ a {\it self-similar Markov} composition structure. 

\par Bernoulli sequences in Example 1 yield \MBR\ composition structures.
Another instance is the renewal sequence
in Example 2, with $Q^\uparrow$ a discrete renewal process.

\par As the terminology is meant to suggest, \MBR\ composition structures are
Markov in the sense of Section \ref{mcs}.
To see this, for each $n$ consider 
a Markov chain $Q_n^\uparrow$ with state-space $[n]\cup\{\infty\}$, such that 
$Q_n^\uparrow$ coincides with $Q^\uparrow$ as long as the latter stays in $[n]$, but jumps to $\infty$ at the time when
$Q^\uparrow$ exits $[n]$. Let $Q_n^\downarrow$ be a time-reversal of $Q_n^\uparrow$, so that $Q_{n,0}^\downarrow$
has the same distribution as the value of $Q_n^\uparrow$ immediately before exiting $[n]$.
The chains $Q_n^\downarrow\,,\,n=1,2,\ldots$ are {\it coherent} in the sense that, for $m\leq n$, when $Q_n^\downarrow$
enters $[m]$ its state has the same distribution as $Q_{m,0}^\downarrow$. Conversely,
if the chains $(Q_n^\downarrow)$ are coherent, their reversals $(Q_n^\uparrow)$ can be organised in a single
`super-chain' $Q^\uparrow$ with state-space $\{1,2,\ldots\}$. So
this property distinguishes the \MBR\ case within the general Markov case.
Another feature characterising the coherent sequence is that there is a common potential function:
for all $n\geq m$ the probability that $Q_n^\downarrow$ visits state $m$ does not depend on $n$.

\par We recall that a {\it stationary regenerative set} \cite{taksar} is the range of a process
$(X+Y_t\,,\,t\geq 0)$ where $(Y_t)$ is a subordinator with L{\'e}vy measure satisfying
\begin{equation}\label{mean}
{\tt m}=\int_0^\infty y\,\nu({\rm d}y)<\infty\,,
\end{equation}
$\drift\geq 0$ is some drift coefficient, $X$ is independent of $(Y_t)$ and has distribution
\begin{equation}\label{ini}
{\mathbb P}(X\in {\rm d}y)={\nu [y,\infty]\over\drift+{\tt m}}\,{\rm d}y+{\drift\over \drift+{\tt m}}\delta_0({\rm d}y).
\end{equation}
Thus, $(X+Y_t)$ is a delayed subordinator with a special choice of distribution for $X$, to make the range
stationary.

\begin{theorem}\label{Mselfsim}
A composition structure ${\cal C}$ is \MBR\ if and only if its associated self-similar
set $Z$ can be presented as $Z=\exp(-W)$ where $W$ is a regenerative set with stationary delay.
\end{theorem}
{\it Proof.}
Follows by combining Theorems \ref{general}
and \ref{mark-gen}.
\endpf

\par We see that self-similarity of Markov composition structures imposes further constraints on the 
decrement matrices
$q$ and $q_*$ in the product formula (\ref{prod}). Thus, $q$ can be associated only with a finite-mean subordinator,
and is given then by (\ref{qPhi}) with $\Phi$ as in (\ref{Phis}) and (\ref{Phinm}).
Similarly, $q_*$ is given by (\ref{q*Psi}) for $\Psi$ as in (\ref{Psinm}), and $A_1$ having distribution
\begin{eqnarray*}
{\mathbb P}(A_1\in {\rm d}x)={\widetilde{\nu} [x,1]{\rm d}x \over (\drift+{\tt m})(1-x)}\,+
{\drift\over \drift+{\tt m}}\delta_0({\rm d}x),
\end{eqnarray*}
where $\widetilde{\nu}$ is the image of $\nu$ under $y\mapsto 1-e^{-y}$. 
By Theorem \ref{A=V} this is also the structural distribution of $Z$, and comparing with (\ref{form}) we
observe that the distribution is of exactly the same type as for the general self-similar $Z$
according to Theorem \ref{strss}.
For the potential function there is a simple formula
\begin{equation}\label{green}
g(j)={1 \over \drift+{\tt m}}{\Phi(j-1)\over j-1}\,,~~{\rm for~~} j>1\,,~~g(1)=1
\end{equation}
which appeared in \cite[p. 86]{sieve} in a special case, and the transition probabilities $f$ of $Q^\uparrow$
are recovered from 
$$q(j-1:j-i)={f(j\,|\,i)g(i)\over g(j)}\,,~~~1\leq i<j\,.$$

\par The relation between a regenerative composition with decrement matrix $q$ and 
the associated \MBR\ composition structure with matrices $q$ and $q_*$ (with $q_*$ given by (\ref{q*Psi}))
is the combinatorial counterpart of the relation between a subordinator and its stationary version.

\subsection{Arrangements} 

A difficult and interesting question is the relation between partition structures 
and their possible {\it arrangements} as
composition structures with certain properties.
Some aspects of this problem were treated in \cite{rcs, rps}.
Although we do not know a simple algorithm
to check 
if the blocks of a given partition may be ordered to produce a \MBR\ composition structure,
we can show the uniqueness.

\begin{proposition} 
\label{unique}
If a partition structure admits an arrangement as a
\MBR\ composition structure, then such arrangement is unique in
distribution.
\end{proposition}
\noindent
The claim follows from the next lemma by recalling 
that the moments of the structural distribution $p(n)$ are determined by the associated 
partition structure.

\begin{lemma} \label{le9} For $({\cal C}_n)$ a \MBR\ composition structure,
for each $n$.
the distribution of ${\cal C}_n$ is uniquely determined 
by the structural moments $p(1)=1,\,p(2)\,,\ldots,\, p(n+1)$.
\end{lemma}
{\it Proof.} A binomial expansion in (\ref{mur1}) shows that $\mu_{n,r}\,, 1\leq r \leq n$, are computable
from $p(1),\ldots,p(n)$. By Theorem \ref{S=L} and the formula
$q_*(n:r)=r\,\mu_{n,r}/n$ also $(q_*(n':r), 1\leq r\leq n'\leq n+1)$ is computable from $p(1),\ldots,p(n+1)$.
Applying (\ref{q*})
we see by induction that the minor of the decrement matrix $(q(n':r), 1\leq r\leq n'\leq n)$
is computable from $p(1),\ldots,p(n+1)$, which taken together with (\ref{prod}) proves the claim.
\endpf 
\noindent
For the regenerative case (when $q=q_*$) we have shown that only the moments $p(2),\ldots,p(n)$ are needed
to recover the distribution of the composition of order $n$
\cite[Proposition 7.1]{rcs}. The explicit formulas are rather involved already in that case.

\par While a \MBR\ arrangement (if any) of a partition structure is unique,
many self-similar composition structures may project onto the same partition structure.
For example, for self-similar $Z$, $M\subset [0,1]$ and $M'$ the reflection of $M$ about $1/2$, both
fragmentation products
$Z\otimes M$ and $Z\otimes M'$ induce the same partition structure, but the composition structures are different,
unless $M\stackrel{d}{=}M'$. 
This implies nonuniqueness in the problem of binary representability of partition structures 
studied in \cite{young}.

\section{The two-parameter family}

We are interested in self-similar composition structures associated with the 
members of the two-parameter family of partition structures \cite{csp}.
For the range of parameters 
$\theta >-\alpha$, $0\leq\alpha<1$ these partition structures 
may be introduced as follows. 

Let $(V_i)$ or $(V_i, i \in I)$ denote a {\em random discrete distribution}, 
that is a collection of random variables indexed by $i$ in some finite or
countably infinite set $I$, with
$$
V_i \ge 0 \mbox{ and } \sum_i V_i = 1 \mbox{ almost surely }.
$$
We use $\{V_i\}$ as an informal notation for multi-set of all non-zero 
values of $V_i$, without regard to how they are indexed by $I$.
Formally, $\{V_i\}$ is encoded by the sequence 
$(\hat{V}_j, j = 1,2, \ldots) := {\rank} (V_i)$ 
meaning that $(\hat{V}_j, j = 1,2, \ldots)$ is the decreasing 
rearrangement of $(V_i)$ with padding by zeros if necessary.
Let us write simply 
\begin{equation}\label{jp0}
\{V_i \} \sim (\alpha,\theta)
\end{equation}
if ${\rank} (V_i)$ has the Poisson-Dirichlet distribution
with two parameters $(\alpha,\theta)$, defined following \cite{pd2, csp}
as the distribution of
${\rank} (\widetilde{V}_i)$ where
\begin{equation}\label{jp1}
\widetilde{V}_1 := W_1 
\end{equation}
has beta$(1-\alpha,\alpha + \theta)$ distribution and for $ i \ge 1$
$$
\widetilde{V}_i := (1-W_1) \cdots (1 - W_{i-1}) W_i
$$
where $W_i$ has beta$(1-\alpha,\alpha + i \theta)$ 
distribution, and the $W_i$ are independent.
It is known 
\cite{pd2} 
that if $\{V_i \} \sim (\alpha,\theta)$ then such
$\widetilde{V}_i$ can be constructed by size-biased random permutation of
$\{V_i \}$. Then $\widetilde{V}_1= V_J$ for a random index $J$ with
$$
\prob (J = j \giv \{V_i \} ) = V_j
$$
while
\begin{equation}\label{jp3}
\{ V_k^\# \}:= \left\{ \frac{ V_i - V_N}{1 - V_N } , ~i \ne J \right\} \sim (\alpha, \alpha + \theta)
\end{equation}
and
$$
\{V_k^\#\} \mbox{ is independent of } \widetilde{V}_1 .
$$
Since $\{V_i \}$ can be measurably recovered 
from
$\widetilde{V}_1$ and
$\{V_k^\#\}$ as 
\begin{equation}\label{jp4}
\{V_i \} = \{ \widetilde{V}_1 \cup \{ 1 - \widetilde{V}_1) V_k^\# \}
\end{equation}
an immediate consequence is
\begin{lemma} {\em \cite[Proposition 35]{pd2}}
\label{pd2}
If $\widetilde{V}_1$ has beta$(1-\alpha,\alpha + \theta)$ distribution and 
$\{V_k^\#\} \sim (\alpha,\theta+\alpha)$, and $\{V_i\}$ is defined by {\rm \re{jp4}}, then
$(V_i) \sim (\alpha,\theta)$ and $\widetilde{V}_1$ is a size-biased pick from
$(V_i)$.
\end{lemma}

\par In \cite{rcs} we established that for $0 \le \alpha<1, \theta \ge 0$
a random discrete distribution $\{V_i\} \sim (\alpha,\theta)$ can be
derived from the interval partition of $[0,1]$ associated with a unique 
regenerative composition structure. Specifically, the image of the L{\'e}vy measure of this $(\alpha,\theta)$ regenerative composition structure
under $y\mapsto 1-e^{-y}$ is
the measure $\widetilde{\nu}$ on $]0,1]$ characterised by
\begin{equation}\label{qnu}
\widetilde{\nu}[x,1]=x^{-\alpha}(1-x)^\theta\,
\end{equation}
and the decrement matrix is
\begin{equation}\label{qnu1}
q(n:r)={n\choose r}{(1-\alpha)_{r-1}\over (\theta+n-r)_r}~{(n-r)\alpha+r\theta\over n}.
\end{equation}
By combining these known results we now obtain the following:
\begin{theorem}\label{2paramet} 
\label{thmalth}
For $0 \le \alpha <1,~\theta>0$ let $Z = \exp(-W)$ where $W$ is the 
stationary version of the regenerative set associated as above with
an $(\alpha,\theta)$ regenerative composition structure.
Then $Z$ is a self-similar Markov random set associated with an 
$(\alpha,\theta-\alpha)$ partition structure. The structural
distribution of $Z$ is beta$(1-\alpha,\theta)$, and $Z$ has the strong
sampling property \rem{ssp}.
\end{theorem}
{\it Proof.}
The structural distribution of $Z$ is read from 
(\ref{qnu}), (\ref{form}) and (\ref{ini}). The construction of $Z$ allows
the application of Lemma \ref{pd2}, 
with $\theta$ replaced by $\theta-\alpha$, to deduce the other conclusions.
\endpf

\par Theorem \ref{thmalth} can also be derived more combinatorially
as follows. Consider the Polya-Eggenberger distributions
$$
q_{\alpha,\theta}(n:r)={n-1\choose r-1}{(\theta+\alpha)_{n-r}(1-\alpha)_{r-1}\over (\theta+1)_{n-1}}\,,~~~~r=1,\ldots,n
$$
and define a function on compositions
\begin{equation}\label{sibi}
\widehat{\pi}_{\alpha,\theta}(\lambda)=\prod_{k=1}^{\ell} q_{\alpha,\theta+(\ell-k)\alpha}(\Lambda_k:\lambda_k)\,,~~~~
\lambda=(\lambda_1,\ldots,\lambda_{\ell})
\end{equation}
where $\Lambda_k=\lambda_1+\ldots+\lambda_k$. 
The formula (\ref{sibi}) is the distribution of the $(\alpha,\theta)$ 
partition structure with parts arranged from right to left in a size-biased order.
The $(\alpha,\theta)$-partition structure is defined then by 
the partition probability function 
obtained by the symmetrisation of the CPF (see \cite{rps} for more details of this procedure):
\begin{equation}\label{pipi}
\pi_{\alpha,\theta} (\lambda^\downarrow)=\sum_{{\rm distinct~}\sigma}\widehat{\pi}_{\alpha,\theta}(\lambda_\sigma) 
\end{equation}
where 
the summation
extends over all
{\it distinct} permutations 
$\lambda_\sigma=(\lambda_{\sigma(1)},\ldots,\lambda_{\sigma(\ell)})$
of parts of composition $\lambda$.
From (\ref{sibi}) and (\ref{pipi}) follows the recursion
\begin{equation}\label{rec-sibi}
\pi_{\alpha,\theta}(\lambda^\downarrow)=\sum_{{\rm distinct}~\lambda_j\in \lambda^\downarrow} 
q_{\alpha,\theta}(n:\lambda_j)\,\,\pi_{\alpha,\theta+\alpha}(\lambda^\downarrow-\lambda_j)
\,,
\end{equation}
where $\lambda^\downarrow$ is a (ranked, unordered) partition of $n$ and 
where $\lambda^\downarrow -\lambda_j$ is the partition $\lambda^\downarrow$ without part $\lambda_j$.
Let ${\cal C}_n$ denote the self-similar Markov composition structure
derived from $S = \exp(-W)$ as in the theorem.
Computing beta integrals to determine the distribution $q_*$ of the last part of ${\cal C}_n$
according to (\ref{q*Psi}) we obtain
\begin{eqnarray*}
q_*(n:r)&=
&{n\choose r}{\mathbb E}\,A_1^r (1-A_1)^{n-r}+{\mathbb E}\,(1-A_1)^n q(n:r) \\
&=& {{n\choose r}\over {\rm B}(1-\alpha,\theta)}
\left({\rm B}(r-1+\alpha,n-r+\theta)+{\rm B}(1-\alpha,n+\theta){(1-\alpha)_{r-1}\over (\theta+n-r)_r}\,\,
{(n-r)\alpha+\alpha\theta\over n}\right)
\end{eqnarray*}
which upon simplification shows that
$q_*=q_{\alpha,\theta-\alpha}$. 
If the last part of ${\cal C}_n$ is $r$ then the rest of ${\cal C}_n$
must be a copy of the $(\alpha,\theta)$
regenerative composition 
of $n-r$,
hence the partition structure can be recovered from

$$
\pi(\lambda^\downarrow)=\sum_{{\rm distinct}~\lambda_j\in \lambda^\downarrow}
q_*(n:\lambda_j)\,\,\pi_{\alpha,\theta}(\lambda^\downarrow-\lambda_j)
$$
which by comparison 
with (\ref{rec-sibi}) shows that $\pi=\pi_{\alpha,\theta-\alpha}$,
in accordance with the conclusion of the theorem.

\begin{corollary}
For $0 \leq \alpha<1,~\theta>-\alpha$ each $(\alpha,\theta)$ partition 
structure has a distributionally unique arrangement as a 
self-similar Markov composition structure. 
\end{corollary}

\par There is an explicit stochastic algorithm which allows, for each $n$, arranging 
an unordered collection of parts of a $(\alpha,\theta)$ partition into a Markovian self-similar composition.
Given a partition $\lambda$ choose a part $\lambda_{j}$ by a size-biased pick and declare it to the right end of the composition under construction. 
Then arrange the rest parts $\lambda-\lambda_j$ one-by-one,
as for the regenerative $(\alpha, \theta+\alpha)$ composition (from right to left), using
the appropriate deletion kernel \cite{rcs, rps}.
Specifically, when the rest partition is $\mu$, the algorithm selects
{\it each} part of size $r$ of $\mu$ with probability
$$ {1\over n}\,\, {(|\mu|-r)\tau +r(1-\tau)\over 1-\tau+(k-1)\tau}$$
where $\tau=\alpha/(2\alpha+\theta)$ and $k$ is the number of parts of $\mu$; 
then the same procedure is applied to the reduced partition, etc.
For example, consider the $(\alpha, 0)$ partition of $n$, assumung it has $\ell$ 
parts, after placing a size-biased pick the rest $\ell-1$ parts should be arranged in a random order, with
all $(\ell-1)!$ orders being equally likely.

\par By the very construction, conditionally
given the parts, the last part is a size-biased pick from all parts: this feature is a combinatorial 
analogue of the 
strong sampling property in Section \ref{block} (as was stated for the $(\alpha,0)$ case
in \cite[Proposition 11 (i)]{bernoulli}).
Algebraically, the combinatorial
strong sampling property amounts to the identity
$$
{q_{\alpha,\theta}(n:\lambda_j)\,\,\pi_{\alpha, \alpha+\theta}(\lambda^\downarrow-\lambda_j)\over \sum_i
q_{\alpha,\theta}(n:\lambda_i)\,\,
\pi_{\alpha, \alpha+\theta}(\lambda^\downarrow-\lambda_i)}=
{\lambda_j\over n}
$$

\subsection{Case $(\alpha=0,~\theta>0)$}

This is case of Example 1, with independent digits and potential function
$$
g(j)={\theta\over j+\theta-1}.
$$
Here $S$ is PPP$(\theta {\rm d}y/y)$. A characteristic feature is that it is the only \MBR\ compostion structure which is right-regenerative.
Indeed, if a random set is both regenerative and stationary regenerative, it is a homogeneous PPP.

\vskip0.5cm

The transition function for the $Q^\uparrow$ chain is
$$
f(j\,|\,i)={\theta (j-2)!\,(\theta)_i\over (i-1)!\,(\theta)_j}
$$

\noindent
{\bf Remark on records.}
The case $\theta=1$ has classical interpretation in terms of indicators of records in a
sequence of i.i.d. random variables with some continuous distribution.
With reference to a question left open in \cite[p. 297]{Bunge},
a similar interpretation exists for any $\theta>0$, but distributions of independent
variables should be different.
One possibility, based on a planar homogeneous Poisson process is the following: 
divide the positive quadrant ${\mathbb R}^2_+$ into vertical strips of widths 
$\beta_j=(\theta)_{j-1}/(j-1)!\,$ and define the variables to be the heights of the 
lowest Poisson atoms in the strips, from left to right.
Elementary algebra shows that, to agree with ESF($\theta$), the collection of $\beta_j$'s 
must be as above up to a common positive factor.
The same distribution of record indicators appears for an independent sample from distributions
$F^{\beta_1},\,\,F^{\beta_2}, \ldots$
where $F$ is an arbitrary continuous distribution on $\mathbb R$.

\subsection{Case $(0<\alpha<1,~\theta=0)$}
The range $S$ of an $\alpha$-stable subordinator 
induces the renewal composition structure of Example 2. This
is the self-similar version
of the regenerative $(\alpha, \alpha)$ composition structure, whose L{\'e}vy measure after the transform $x = 1-e^{-y}$ 
is $\widetilde{\nu}_{\alpha,\alpha}$ defined by
$$\widetilde{\nu}_{\alpha,\alpha}[x,1]=x^{-\alpha}(1-x)^\alpha \,,$$
hence the potential function is

$$
g(j)={\Phi_{\alpha,\alpha}(j-1)\over {\tt m}_{\alpha,\alpha}\,(j-1)}={(\alpha)_{j-1}\over (j-1)!}.
$$
(where ${\tt m}_{\alpha,\alpha}$ is the mean value as in (\ref{mean})).

\par The induced composition structure is self-similar Markov
as well as left-regenerative.
Thus $S\cap[0,1] = 1 - e^{-R}$ for
$R$ the range of another, killed, subordinator, as detailed in \cite{rcs}.
The combination of the two regeneration properties is characteristic:

\begin{proposition}
If a composition structure $({\cal C}_n)$ is both Markov self-similar
and left regenerative, then $({\cal C}_n)$ is the $(\alpha,0)$
composition structure derived by sampling from the range of some
$\alpha$-stable subordinator.
\end{proposition}
{\it Proof.} Let $Z$ be the set in Kingman's representation of $({\cal C}_n)$.
The left regeneration property implies that $Z$ is the range of a multiplicative subordinator 
$1-e^{-A}$, where $A$ is some subordinator. 
On the other hand, by Theorem \ref{Mselfsim}, $Z=e^{-B}$ for $B$ some stationary delayed subordinator,
hence $Z$ has a nontrivial meander with positive probability, which implies that $A$ has a positive killing rate.
Let $Z_0$ be the set $Z$ conditioned on zero meander, which is the range of the multiplicative 
subordinator $1-e^{-A_0}$, for $A_0$ the version of $A$ without killing. 
Then, of course, $Z_0=e^{-B_0}$ for $B_0$ the version of $B$ but with zero delay.
It follows that the composition structure induced by $Z_0$ is both left- and right-regenerative,
that is both sets $Z_0$ and $1-Z_0$ are multiplicatively regenerative.
By \cite[Theorem 12.1 and Corollary 12.2]{rcs}, $Z_0\ed 1-Z_0$, 
$Z_0$ is
the zero set of a Bessel bridge,
and the composition structure induced by $Z_0$ is of type $(\alpha,\alpha)$.
By Theorem \ref{2paramet} the stationary version of this composition structure is of type $(\alpha,0)$,
and $Z$ is the range (restricted to $[0,1]$) of some $\alpha$-stable subordinator.
\endpf

\noindent
{\bf Remark.} This result complements the characterisation 
of $(\alpha,\alpha)$ regenerative composition structures in \cite[Theorem 12.1]{rcs}.
Apparently,
the assumption of the Markov property can be omitted,
{\it i.e.} it seems 
sufficient to assume only that $({\cal C}_n)$ is right-consistent.
That the Markov property follows is not obvious, because
the left-regeneration property of $({\cal C}_n)$ does not imply the right Markov property of the composition
in the sense of Definition \ref{Mcs} (which requires time-homogeneity of the Markov chain).
Still, a plausible argument is the following.
As above, define left-regenerative (multiplicatively) $Z_0$ by conditioning on zero meander
(a limiting procedure required to justify this definition is obvious).
Fix $x\in ]0,1[$ and condition on $x\in Z_0$, then, because $Z$ is self-similar, $[0,x]\cap Z_0\ed x\,Z_0$.
But by the left-regeneration (multiplicative) property, $[0,x]\cap Z_0$ is independent of $[x,1]\cap Z_0$,
whence the right-regeneration (multiplicative) property. Then the conclusion is above.
A loose point in this argument is the conditioning on the zero event $x\in Z_0$

\subsection{Case $(\alpha,\alpha)$}
For this partition there is a regenerative arrangement
(the composition structure induced by the Bessel bridge) 
and another \MBR\ arrangement. The latter 
is the self-similar version of the regenerative $(\alpha,2\alpha)$ composition.

\subsection{General $0<\alpha<1, \,\theta>-\alpha$}

Explicit construction of the self-similar Markov composition structure associated with the $(\alpha, \theta)$ 
partition structure exploits the fragmentation product introduced in Section \ref{fp}.
One ingredient is the Poisson process
$Z$ PPP($\theta{\rm d}y/y$), $\theta>0$, restricted to $[0,1]$.
Another factor is the set 
$M'=1-M\cap [0,1]$ 
obtained by reflecting the range $M$ of the $\alpha$-stable subordinator.
The self-similar set is defined then as the fragmentation product $Z\otimes M'$, 
and the induced composition
is the \MBR\ version of partition $(\alpha,\theta -\alpha)$.
Conditioning on zero meander will produce a set corresponding to $(\alpha,\theta)$ regenerative composition,
as in \cite{rcs}.

\par Unlike $M$, the set $M'$ exploited here has the leftmost meander interval.
The fragmentation product $Z\otimes M$ was introduced in \cite{RDD}; 
the resulting composition structure is right-consistent but not Markovian.

\def\cprime{$'$} \def\polhk#1{\setbox0=\hbox{#1}{\ooalign{\hidewidth
\lower1.5ex\hbox{`}\hidewidth\crcr\unhbox0}}} \def\cprime{$'$}
\def\cprime{$'$} \def\cprime{$'$}
\def\polhk#1{\setbox0=\hbox{#1}{\ooalign{\hidewidth
\lower1.5ex\hbox{`}\hidewidth\crcr\unhbox0}}} \def\cprime{$'$}
\def\cprime{$'$} \def\polhk#1{\setbox0=\hbox{#1}{\ooalign{\hidewidth
\lower1.5ex\hbox{`}\hidewidth\crcr\unhbox0}}} \def\cprime{$'$}
\def\cprime{$'$} \def\cydot{\leavevmode\raise.4ex\hbox{.}} \def\cprime{$'$}
\def\cprime{$'$} \def\cprime{$'$} \def\cprime{$'$}

\end {document}